\documentclass[12pt,amsmath,amssym,latexsymb,leqno]{article}

\usepackage{amssymb,latexsym}

\usepackage{graphicx}

\usepackage{psfrag}

\usepackage[final,dvips]{epsfig}

\usepackage{verbatim}

\usepackage{amsfonts}

\def\newtext{}

\makeatletter
\def\@begintheorem#1#2{\trivlist \item[\hskip \labelsep{\sc #1\ #2.}]\it}
\def\@opargbegintheorem#1#2#3{\trivlist
      \item[\hskip \labelsep{\sc #1\ #2.\ (#3)}]\it}
\def\@sect#1#2#3#4#5#6[#7]#8{\ifnum #2>\c@secnumdepth
     \let\@svsec\@empty\else
     \refstepcounter{#1}\edef\@svsec{\csname the#1\endcsname.\hskip 1em}\fi
     \@tempskipa #5\relax
      \ifdim \@tempskipa>\z@
        \begingroup #6\relax
          \@hangfrom{\hskip #3\relax\@svsec}{\interlinepenalty \@M #8\par}%
        \endgroup
       \csname #1mark\endcsname{#7}\addcontentsline
         {toc}{#1}{\ifnum #2>\c@secnumdepth \else
                      \protect\numberline{\csname the#1\endcsname}\fi
                    #7}\else
        \def\@svsechd{#6\hskip #3\relax  
                   \@svsec #8\csname #1mark\endcsname
                      {#7}\addcontentsline
                           {toc}{#1}{\ifnum #2>\c@secnumdepth \else
                             \protect\numberline{\csname the#1\endcsname}\fi
                       #7}}\fi
     \@xsect{#5}}
\def\section{\@startsection {section}{1}{\z@}{-1.5ex plus-1ex minus
    -.2ex}{-2.5ex plus.2ex}{\reset@font\bf}}
\def\subsection{\@startsection{subsection}{2}{\z@}{-3.25ex plus-1ex
    minus-.2ex}{-1.5ex plus.2ex}{\reset@font\sl}}

\makeatother

\newtheorem{predf}{Definition}[section]

\newtheorem{prex}[subsection]{Example}

\newtheorem{theo}[subsection]{Theorem}

\newtheorem{prop}[subsection]{Proposition}

\newtheorem{coro}[subsection]{Corollary}

\begin{document}

\renewcommand{\a}{{\alpha}}
\renewcommand{\b}{{\beta}}
\newcommand{\g}{{\delta}}
\renewcommand{\d}{{\gamma}}

\newcommand{\C}{{\Bbb C}}
\newcommand{\R}{{\Bbb R}}
\newcommand{\D}{{\Bbb D}}
\newcommand{\re}{{\rm Re \,}}
\newcommand{\im}{{\rm Im \,}}
\newcommand{\e}{\varepsilon}
\newcommand{\voltp}[1]{V^{(#1)}}
\newcommand{\currp}[1]{I^{(#1)}}
\newcommand{\forwp}[1]{U^{(#1)}}
\newcommand{\Rn}[1]{{\Bbb R}^{#1}}
\def\class{L_{\infty,\alpha}}

\def \beq {\begin {eqnarray}}
\def \eeq {\end {eqnarray}}
\def \ba {\begin {eqnarray*}}
\def \ea {\end  {eqnarray*}}
\def \p {\partial}
\def \tilde {\widetilde}
\def \hat {\widehat}
     \newcommand{\be}{\begin{equation}}   
\newcommand{\ee}{\end{equation}}
\newcommand{\ben}{\[}                
\newcommand{\een}{\]}
\newcommand{\bea}{\begin{eqnarray}}  
\newcommand{\eea}{\end{eqnarray}}
\newcommand{\bean}{\begin{eqnarray*}}
\newcommand{\eean}{\end{eqnarray*}}
\parindent=.5cm

\newcommand{\F}{{\mathcal{F}}}
\newcommand{\eps}{\varepsilon}
\newcommand{\nbr}{\mathcal{N}}
\newcommand{\tr}{^{\rm T}}
\newcommand{\ith}{^{\rm th}}
\newcommand{\domain}{\Omega}
\newcommand{\compdomain}{\Omega_m}
\newcommand{\modeldomain}{\tilde \Omega}

\def\mtrx#1#2{
  \left(
    \begin{array}{#1}
      #2
    \end{array}
  \right)}

\title{The inverse conductivity problem with an imperfectly known boundary in three dimensions}

\author{
Ville Kolehmainen\footnote{Department of physics, University of Kuopio, 
P.O.Box 1627, 70211 Finland},\quad Matti Lassas\footnote{Department of mathematics, 02015 Helsinki University of Technology, Finland},\quad 
Petri Ola\footnote{Department of mathematics, 00014 University of Helsinki, Finland}
\\
}
\date{}
\maketitle

{\bf Abstract.} We consider the inverse conductivity problem in
a strictly convex domain whose boundary is not known. 
Usually the numerical reconstruction from the measured current and voltage data is done assuming the 
domain has a known fixed geometry. However, in practical applications the geometry of the domain is usually not known. This introduces an error, and effectively changes the problem into an anisotropic one. The main result of this paper is a uniqueness result  characterizing the isotropic conductivities on convex domains in terms of measurements done on a different domain, which we call the model domain, up to an affine isometry. 
As data for the inverse problem, we assume 
 the Robin-to-Neumann map and the contact impedance function 
on the boundary of the model domain to be given.
 Also, we present a minimization algorithm based on the use
of Cotton--York tensor, that finds the pushforward of the isotropic conductivity to our model domain, and also finds the boundary of the original domain up to an affine isometry. This algorithm works also in dimensions higher than three, but then the Cotton--York tensor has to replaced with the Weyl--tensor.
\medskip

{\bf AMS classification:} 35J25, 35R30, 58J32.
\medskip

{\bf Keywords:} Inverse conductivity problem, electrical impedance
tomography, unknown boundary, Cotton--York tensor.
\bigskip

\section{Introduction.} 

We consider the electrical impedance tomography problem (EIT for short), 
i.e.\ the determination of the unknown isotropic conductivity distribution inside a domain in $\R ^3$, for example the human thorax, from voltage and current measurements made on the boundary. Mathematically this is 
formulated as follows: Let $\Omega$ be the measurement domain, and 
denote by $\gamma$ the bounded and strictly positive function describing the conductivity in $\Omega$. 
The voltage potential $u$ satisfies in $\Omega$ the 
equation
\begin{equation}\label{johty}
\nabla\cdot \gamma\nabla u = 0.
\end{equation}
To uniquely fix the solution $u$ it is enough to give its value on the boundary. Let this be $f$.
{\newtext In the idealized case, when the contact impedance
of the measurement device is zero, one measures  for
all voltage distributions $u|_{\p M}=f$ on the boundary the corresponding
current flux through the boundary, $\gamma \partial y/\partial \nu$,
 where $\nu$ is the exterior unit normal to 
$\p\Omega$.} Mathematically this amounts to the knowledge of the Dirichlet--Neumann map $\Lambda$ 
corresponding to $\gamma$, i.e. the map taking the Dirichlet boundary values to the 
corresponding Neumann boundary values 
of the solution to (\ref{johty}),
\[
\Lambda: \ \ u|_{\p M}\mapsto \gamma \frac{\partial u}{\partial \nu}
\]
{\newtext The Calder\'on's inverse problem is then to to reconstruct $\gamma$ from
 $\Lambda$.
The problem was originally proposed by Calder\'on \cite{ca} in 1980 
and} then solved in dimensions three and higher
for isotropic conductivities which are $C^\infty$--smooth in \cite{SU}
{\newtext and \cite{Nachman}}. 
The smoothness requirements
have been since relaxed, and currently the best known result is \cite{PPU} with unique determination of conductivities in $W^{3/2,\infty}$, 
see also \cite{GLU} for a somewhat different approach to the lack of smoothness.  In two dimensions the first global result is due to Nachman ({\cite{na2}), and later Astala and
P\"aiv\"arinta showed in \cite{ap} that uniqueness holds also for general isotropic $L^\infty$--conductivities. For the corresponding anisotropic case, see
\cite{alp,LU,LUT,LeeUhlmann} and numerical implementations of the methods with
simulated and real data, see \cite{Siltanen3,Siltanen2,Siltanen1}.
\smallskip

{\newtext Assuming that 
the measured Dirichlet-to-Neumann map $\Lambda_{\rm{meas}}$ is given,
an  often used method to solve the  EIT--problem is to minimize}
\[
\|\Lambda_{\rm{meas}} - \Lambda _{\sigma}\|^2 +\alpha \|\sigma\|_X ^2
\]
for $\sigma$ defined in terms of some triangulation of $\Omega$ and $\|\cdot\|_X$ is some regularization norm; here $\Lambda _{\sigma}$ is the Dirichlet--Neumann map corresponding to the conductivity $\sigma$.
One then also fixes the geometry of $\Omega$ by assuming that it is for example a ball or an ellipsoid.
Now, if our measurements have no error a Bayesian interpretation
of this problem as a search of an MAP-estimate suggests that $\alpha=0$.
{\newtext Usually, the given data $\Lambda_{\rm{meas}}$ 
does not correspond to
any isotropic conductivity in the model domain. The reason for this is
that there is no conformal map deforming the original domain 
to the model domain. Therefore, in solving the minimization 
problem} we obtain an incorrect solution
$\sigma$. This means that a systematic error in modeling
causes a systematic error to the reconstruction.
In particular, if we consider linearization 
$\gamma=\gamma_0+\e\gamma_1$ where
$\gamma_0$ is given known background conductivity and $\e$ is small,
it seems that a  localized perturbation $\gamma_1$ 
gives a  reconstruction $\sigma=\gamma_0+\e\sigma_1$ 
where the  reconstructed perturbation $\sigma_1$ is not localized. 
This is clearly seen in brain-activity measurements, see \cite{gersing96} and \cite{kolehmainen97e}.
\smallskip

This work is continuation of \cite{Kolehmainen} where the corresponding question in two dimensions was studied: we proved that on the model domain there is a unique (anisotropic) conductivity with minimal anisotropy. This follows from a result of Strebel saying that among all quasiconformal self-maps of the unit disk with a fixed boundary value there is a unique one with minimal complex dilation. In higher dimensions there are several new issues. First of all, the 
non-uniqueness  due to anisotropy is not understood, except 
{\newtext in the case when both the domain and the conductivity function are
 the real analytic (\cite{LeeUhlmann}, \cite{Lionheart})}. Also, as we already mentioned, in the plane case one could use the theory of quasiconformal maps to break the non-uniqueness. The higher dimensional analogue of this is unknown. Finally, there is no analogue of the Riemann mapping theorem that we could use.
\smallskip

The structure of this paper is the following. In the first part, consisting of sections 2--4,  we present the uniqueness results we have on the problem. It is worth noting that we choose to work with the Robin--to--Neumann (RN) map instead of the Dirichlet--to--Neumann (DN) map described above. Mathematically they are equivalent, as we will show, but the RN--map is a better model for the actual measurement configuration since it takes into account the contact impedances at $\partial \Omega$
\cite{Somersalo92}. 
{\newtext Also, we assume the the function modeling the contact impedances
of the electrodes is known.}
There are two key ideas how we compensate for our lack of understanding of the full anisotropic problem. First thing is to note that if an isotropic conductivity is pushed forward by a diffeomorphism, the resulting conductivity is still conformally flat, and in three dimensions this is equivalent with the vanishing of the Cotton--York tensor. Secondly, we assume that our original domain is strictly convex, and then the Cohn--Vossen theorem can be used to determine the original boundary $\p\Omega$ up to rigid motions.
\smallskip

In the second part we develop an algorithm for finding the 
shape of the domain $\Omega$ and the conductivity inside
using a minimization technique. Important feature is
that we do not have to construct an embedding of the
boundary to the Euclidean space. We plan to report on the numerical implementation of 
our algorithm in a separate article.
\bigskip

\section{Measurements.} \label{theorysec}

Let $\Omega\subset \R^n$, $n\geq 3$ be a strictly convex domain, and
denote by $\gamma = (\gamma ^{ij}(x))_{i,j=1}^n$ the symmetric real valued
matrix describing the conductivity in $\Omega$.
We assume that the matrix is bounded from above
and from below, that is, for some $C,c>0$ we have
\beq\label{ehto 1}
c\|\xi \|^2\leq \langle\xi,\gamma (x)\xi\rangle \leq C\|\xi \|^2,\quad
\hbox{ for all }x\in \Omega.
\eeq
We will state the precise smoothness of $\gamma$ later. 
{\newtext We start by consider the 
EIT problem with continuous boundary data. Instead of the
Dirichlet-to-Neumann map we will use the Robin-to-Neumann map
defined below that corresponds better to the measurements
done in practice. We discuss later in this
section the relation of the continuous model and the 
electrode measurements made in practice.}

For the electrical potential $u$ we write the model
\beq\label{joht}
& &\nabla\cdot \gamma\nabla u = 0,\quad\quad  x \in \Omega,  \\
& &(z\nu\cdotp \gamma \nabla u + u)|_{\p \Omega}=h,
\eeq
where $h$ is the Robin-boundary value of the potential
and $z$ is a function describing the contact impedance on the boundary. The contact impedance models the impedance that is caused by 
electro-chemical phenomena at the interface of the skin and the 
measurement electrodes in practical measurements \cite{chengea89}. 
\smallskip

In mathematical terms, the perfect boundary measurements are modeled by
the Robin-to-Neumann map $R=R_{z,\gamma}$ given by
\ba
R: h \mapsto \nu\cdotp \gamma \nabla u|_{\p \Omega}
\ea
that maps the potential on the boundary
to the current across the boundary. 
{\newtext Next we relate this continuous model
for measurements done in practice.}
\smallskip

The physically realistic measurements are usually modeled by}
the following {\em complete electrode model}  
(see \cite{chengea89, Somersalo92}): Let $e_j\subset \p \Omega$,
$j=1,\dots,J$
be disjoint open sets of the boundary modeling the electrodes that are used
for the measurements. Let $u$ solve the equation
\beq
\label{poissoneq1}
& &\nabla\cdotp \gamma \nabla v=0\quad\quad\hbox{in }\Omega,\\
\label{poissoneq2}
& &z_j\nu\cdotp \gamma \nabla v+v|_{e_j}=V_j,\\
\label{poissoneq3}
& &\nu\cdotp \gamma \nabla v|_{\p \Omega\setminus \cup_{j=1}^J e_j}=0,
\eeq
where $V_j$ are constants
 representing electric potentials on electrode $e_j$. 
Then, one measures the currents observed
on the electrodes, given by
\ba
I_j=\frac 1 {|e_j|}\int_{e_j}\nu\cdotp \gamma \nabla v(x)\,ds(x),\quad j=1,\dots,J.
\ea
Thus the electrode measurements are given by map
$E:\R^J\to \R^J$, $E(V_1,\dots,V_J)=(I_1,\dots,I_J)$.
We say that $E$ is the electrode measurement matrix
for $(\p \Omega,\gamma,e_1,\dots,e_J,$ $z_1,\dots,z_J)$.
\smallskip

The complete electrode model can alternatively be defined as follows: the Robin-to-Neumann map $R_\eta$ is given by
$
R_\eta f=\nu\cdotp \gamma \nabla u|_{\p \Omega}
$
where  $u$ is the solution of
\beq \label{N 1}
& &\nabla\cdotp \gamma \nabla u=0\quad\quad\hbox{in }\Omega\\
& &z\nu\cdotp \gamma \nabla v+\eta v|_{\p \Omega}=h, \nonumber
\eeq
where $z\in C^\infty(\p \Omega)$ is such that
its restriction to the electrode $e _j$ is equal to the constant $z_j$ and $\eta=\sum_{j=1}^J \chi_{e_j}$,
where $\chi_{e_j}$ is the characteristic function of electrode $e_j$. 
\smallskip

We associate to the electrode measurement matrix and to the complete electrode model also the corresponding quadratic forms
$E:\R^J\times
\R^J\to \R$ and $R_\eta :H^{-1/2}(\p \Omega)\times H^{-1/2}(\p \Omega)\to \R$
given by
\beq\label{forms}
E[V,\tilde V]=\sum_{j=1}^J (EV)_j\tilde V_j |e_j| ,\quad
R_\eta [h,\tilde h]=\int_{\p \Omega} (R_\eta h)\,\tilde h\,ds.
\eeq
These have the following simple relation to each other:
Let $S=\hbox{span}(\chi_{e_j}:\ j=1,\dots,J)\subset H^{-1/2}(\p \Omega)$
and define $M:V=(V_j)_{j=1}^J\mapsto \sum_{j=1}^J V_j\chi_{e_j}$
to be a map $M:\R^J\to S$. Then
\beq\label{discretization}
E[V,\tilde V]=R_\eta [MV,M\tilde V].
\eeq
{\newtext By (\ref{discretization}), the electrode measurement matrix
can be viewed as the discretization of the form $R_\eta$. By increasing the
number of the electodes and making the gaps between them smaller, we can assume that $\eta\to 1$. In this case $R_\eta$ approximates the Robin-to-Neumann map
$R_{\gamma,z}$. 
Note that $E(V,V)$ corresponds to the power needed to maintain the voltages
$V$ in electrodes.
}
\smallskip

In practical EIT experiments, one places a set of measurement electrodes on the boundary $\partial \Omega$,
e.g., around the chest of the patient. All the traditional approaches to the numerical EIT reconstruction
assume that the shape of the domain $\Omega$ is known and the only unknown is the conductivity $\gamma$.  
However, in most EIT experiments the boundary of the body $\Omega$ is not known accurately and
since there are no practically reliable measurement 
methods available for the determination of the boundary,
the EIT image reconstruction problem is typically solved
using an approximate model domain $\tilde\Omega$,
which represents our best guess for the shape of the true body $\Omega$.
However, it has been noticed that the use of slightly incorrect
model for the body $\Omega$ in the numerical reconstruction
can lead to serious artefacts in reconstructed images 
\cite{kolehmainen97e,adler96b,gersing96}. This situation is our paradigm for the
EIT problem when the boundary is unknown. {\newtext Next
we analyze how the deformation of the domain affects measurements.}

\bigskip

\section {Deformations of the domain.}
In this section we analyze the behavior of the electrode models under a diffeomorphism.
Let's consider first the Robin--to--Neumann map $R$. The corresponding quadratic form, which we still denote by $R$, is given on the diagonal by
\beq\label{form B}
\quad\quad\quad R[h,h]=\int_{\p \Omega}
(u+z\nu\cdotp \gamma \nabla u)\nu\cdotp \gamma \nabla u\,dS_E
=\int_{ \Omega} \gamma \nabla u\cdotp \nabla u\,dx+
\int_{\p \Omega}z\,
|\nu\cdotp \gamma \nabla u|^2\,dS_E
\eeq
where $h\in H^{-1/2}(\p \Omega)$, $u$ solves (\ref{N 1}),
and $dS_E$ is the Euclidean volume form (or area) of $\p \Omega$. 
The value $R[h,h]$
corresponds to the power needed to maintain the current $h$ on the boundary.
From the mathematical viewpoint,
using the (incorrect) model domain $\tilde\Omega$ instead of the original
domain $\Omega$ can be viewed as a deformation of the original domain. Thus,
let us next consider what happens to the conductivity equation
when the domain $\Omega$ is deformed to $\tilde \Omega$.
Assume that $F:\Omega\to \tilde \Omega$ is a
sufficiently smooth orientation preserving map with
sufficiently smooth inverse $F^{-1}:\tilde \Omega\to \Omega$.
Let $f:\p \Omega\to \p \tilde \Omega$ be the restriction of $F$ on the
boundary. When  $u$ is a solution of $\nabla\cdotp \gamma\nabla u=0$ in $\Omega$,
then $\tilde u(\tilde x)=u(F^{-1}(\tilde x))$ and
$\tilde h(x)=h(f^{-1}(x))$
satisfy the conductivity
equation
\beq
& &\nabla\cdotp \tilde \gamma\nabla \tilde u=0,\quad \hbox{in }\tilde
\Omega,\\
& &\tilde z\tilde \nu\cdotp \tilde \gamma \nabla \tilde u+\tilde u|_{\p \tilde \Omega}=\tilde
h. \nonumber
\eeq
Here 
$\tilde h(x)=h(f^{-1}(x))$, $\tilde \nu$ is
the unit normal vector of $\p \tilde\Omega$,  $\tilde z$ is the deformed contact impedance
 and
$\tilde \gamma$ is the conductivity
\beq\label{pushforward}
\tilde \gamma (x) = \left.\frac{F'(y) \, \gamma(y)\, (F'(y))^{\rm T}}{|\det
F'(y)|}
\right|_{y= F^{-1}(x)},
\eeq
where $F'=DF$ is the Jacobian of the map $F$.
Note that even if $\gamma$ is isotropic, i.e., scalar valued the deformed
conductivity $\tilde \gamma $ can be anisotropic  i.e., matrix valued.
\smallskip

To determine the deformed contact impedance
 $\tilde z$, we consider the corresponding invariant $(n-1)$-form
\ba
J
:=\nu\cdotp \gamma \nabla u\,dS_E\in \Omega^{n-1}(\p \Omega)
\ea
corresponding 
to the current flux through the boundary. 
Next we denote $\tilde x=F(x)$. 
A straightforward application of chain rule gives that
\[
\tilde\nu \cdot\tilde \gamma\nabla \tilde u|_{\p \tilde \Omega} = 
\left((\det DF)^{ -1}\nu\cdot \nabla u\right) \circ f^{-1}|_{\p \tilde \Omega}
\]
since $F$ was orientation preserving and $DF$ 
is the Jacobian of $F$ in
{\newtext boundary normal coordinates
associated to the surface $\p \Omega\subset \R^n$.
In these coordinates $\det DF|_{\p \Omega}=
\det Df$, where $\det Df$ is the determinant of the differential
of the the boundary map $f:\p \Omega\to \p \tilde \Omega$.
We note that $(\det Df\circ f^{-1})f_*(dS_E)=d\tilde S_E$,
where $dS_E$ and $d\tilde S_E$ are Euclidean volume forms of
$\p \Omega$ and $\p \tilde \Omega$, respectively. }
Hence, $z\nu\cdot \nabla u$ transforms as an invariantly defined function when the contact impedance is interpreted as a density, i.e. 
\beq\label{eq: transformation of z}
\tilde z(\tilde x) =  (\det Df(x))z(x)
\eeq
where $f(x) = \tilde x$.
Now we see that the boundary measurements are invariant: When
$f:\p \Omega\to \p\tilde \Omega$ is the restriction of $F:\Omega\to 
\tilde \Omega$,
we say that the map $\tilde R = f_*R_{z,\gamma}$, defined by
\[
((f_*R_{z,\gamma})h)(x)= \left.(R_{z,\gamma}(h\circ f))(y)
\right|_{y= f^{-1}(x)},\quad h\in H^{1/2}(\p \tilde\Omega)
\]
has that property that $\tilde R=R_{\tilde z,\tilde \gamma}$.
We call $\tilde R $ the
push forward of $R_{z,\gamma}$ by $f$. 
\smallskip

Is is also worth noting that in formula (\ref{form B}) the integral over $\Omega$, as well as the integral
over the boundary are  invariant because of the deformation rule
(\ref{eq: transformation of z}) for the contact impedance $z$, that is, we have
\ba
R[h,h']=\tilde R[h\circ f^{-1}, h'\circ f^{-1}],
\ea
for $h,h'\in H^{-1/2}(\p \Omega)$.
\bigskip

\section{Uniqueness results}

Now we are ready to give the 
exact
set--up of the problem we consider:
We want to recover an image of
the unknown conductivity $\gamma$
in $\Omega$ from the measurements of Robin-to-Neumann map,
and we  assume {\em a priori} that $\gamma$ is isotropic.
We assume $z$, $\partial \Omega$ and $R$ are not known.
Instead, let $\modeldomain$, called the model domain,
be our best guess for the domain and let
$f_m:\partial \Omega \to\partial \Omega_{m}$
be a diffeomorphism modeling the approximate knowledge of the boundary.
\smallskip

As the data for the inverse problem, we assume that we
are given the boundary of the model domain
$\p \tilde\Omega$, the function $z\circ f^{-1}$ corresponding
to the contact impedance of electrodes, and
the Robin-to-Neumann map $\tilde R=(f_m)_* R$.
{\newtext Note the discrete analog of this data is to know
the voltage-to-power form $V\mapsto E(V,V)$ and the contact impedances
of the electrodes, but not the location of the electrodes or 
the boundary of the domain.
It is reasonable to assume  that 
the contact impedance $z\circ f_m^{-1}$ on the boundary of the model domain is 
known  since 
we can observe and set up the contact impedances of the electrodes
 the way we want. } 
Hence we have on the boundary of our
model domain $\p \tilde\Omega$ a boundary map $\tilde R$ that
does not generally correspond to any isotropic conductivity.
Furthermore, we saw above that there are many anisotropic
conductivities for which Robin-to-Neumann map is the given map 
$\tilde R$. Next we show that the existence of the ``underlying``
isotropic conductivity in $\Omega$ gives the uniqueness in $\tilde \Omega$
up to a diffeomorphism and that the domain $\Omega$ and the isotropic
conductivity on it can be uniquely determined.
\smallskip

\begin{theo}\label{seainoaoikea}
Let $\Omega\subset \R^n$, $n\geq 3$ be a bounded, strictly convex, 
$C^{\infty}$--domain. Assume 
that $\gamma \in C^{\infty}(\overline
\Omega)$ is an isotropic conductivity, $z\in C^\infty(\p \Omega)$,
$z>0$ be a contact impedance, and $R_{\gamma,z}$ the
corresponding Robin-to-Neumann map. 
Let
$\tilde\Omega$ be a model of the domain satisfying the same regularity 
assumptions as $\Omega$, and 
$f_m :\p\Omega\to \p\tilde\Omega$ be  a $C^{\infty}$--smooth orientation preserving diffeomorphism.

Assume that we are given $\p\tilde \Omega$, the values 
of the contact impedance $z(f_m^{-1}(\tilde x))$,
$\tilde x\in \p\tilde\Omega$  and the map $\tilde R=(f_m)_*R_{\gamma,z}$.
Then we can determine $\Omega$ up to a rigid motion $T$ and the 
conductivity $\gamma\circ T^{-1}$ on the reconstructed domain $T(\Omega)$.
\end{theo}
\smallskip

\noindent We recall also that
rigid motion is an affine isometry $T:\R^n\to \R^n$.
\smallskip

{\bf Proof.}
Assume we are given $\tilde R$ and the values of the contact
impedance, that is, the function
$z(f_m^{-1}(\tilde x))$.
{\newtext Let  $F_m:\Omega\to \tilde \Omega$ be 
an orientation preserving diffeomorphism satisfying
$F_m|_{\p \Omega}=f_m$.} 
As noted before, $\tilde R=R_{\tilde z,\tilde \gamma}$
where $\tilde z(x)=\det(Df_m)z(f_m^{-1}(x))$ 
is the contact impedance on $\p \modeldomain$
and $\tilde \gamma=(F_m)_*\gamma$ is the push forward of $\gamma$ in $F_m$.
The Robin--to--Neumann map is a classical pseudodifferential operator of order zero, with principal symbol $1/\tilde z$, and hence 
$\tilde R$ determines $\tilde z$. Since we assume also $z\circ f^{-1}$ known, we can determine the determinant $\det Df_m$; note that this gives the change of boundary area under deformation $f_m$. For the rest of the proof denote $\beta = \det Df_m$.
Also, this implies that we can find the Dirichlet-to-Neumann  {\newtext 
map $\Lambda_{\tilde \gamma}=(\tilde R^{-1}-\tilde zI)^{-1}$ on $\p \tilde \Omega$},
that is, the map taking the Dirichlet boundary values to Neumann boundary values.
\smallskip

The Riemannian metric corresponding to
the isotropic conductivity $\gamma=\gamma(x) I$ in $\Omega$ is given by 
\ba
g_{ij}(x)=\det(\gamma(x) I)^{1/(2-n)}\,(\gamma(x) I)^{-1}=\gamma(x)^{2/(n-2)}\delta_{ij}.
\ea  
Then, if $\Delta _g$ is the Laplace--Beltrami--operator corresponding to the metric $g$, we have
$\Delta _g = |g|^{-1/2} \nabla\cdot \gamma\nabla$m
where $|g|=\det(g_{ij})$.
This metric is an invariant object and
in the deformation $F_m$ it is transformed to the metric
$\tilde g=(F_m)_*g$ in $\tilde \Omega$. By \cite{LeeUhlmann},
the Dirichlet-to-Neumann map $\Lambda_{\tilde \gamma}$ 
determines the metric tensor
$\tilde g_{ij}$ on the boundary in the boundary normal coordinates.
In particular, if we consider $\p \tilde \Omega$ as a submanifold
of $\R^n$ with the metric $\tilde h=\tilde i^*(\tilde g)$ inherited from $(\tilde 
\Omega,\tilde g)$ where $\tilde i:\p \Omega\to \Omega$ is the identity map,
we see that our boundary data determines the metric $\tilde h$ on 
$\p \tilde \Omega$. Let now metric $h=i^*(g)$ be the corresponding
metric on $\p \Omega$, where $i:\p \Omega\to \Omega$ is the identical
embedding. Then we have
\beq\label{eq: push h}
\tilde h=(f_m)_*h,\quad h=\gamma^{2/(n-2)}h^E 
\eeq
where $h^E$ is the Euclidean metric of $\p \Omega$.
Denote by $\tilde h^E=(f_m)_*h^E$ the metric tensor on $\p\tilde\Omega$, i.e. the push--forward of the Euclidean metric of $\p \Omega$ by $f_m$.
{\newtext Recall that
 $d S_E$ and  $d\tilde S_E$ are the Euclidean volume forms
of $\p \Omega$ and $\p \tilde \Omega$, respectively.}
Then
the Riemannian volume forms $dS_{\tilde h}$ and $dS_{h}$ of the metrics $\tilde h$
and $h$ respectively satisfy
\ba
dS_{\tilde h}=(f_m)_*(dS_{h})=\gamma(f_m^{-1}(\tilde x))(f_m)_*(dS_E)=
(\gamma\beta)\circ f_m^{-1}(\tilde x)\, d\tilde S_E
\ea
on $\tilde \p \Omega$. 
As $\beta$ was already determined,
this shows that we can find $\gamma(f_m^{-1}(\tilde x))$,
$\tilde x\in \p\tilde \Omega$ and hence by (\ref{eq: push h})
we can determine the metric 
\ba
\tilde h^E= \gamma(f_m^{-1}(\tilde x))^{-2/(n-2)}\tilde h. 
\ea
In other words, if we consider $\p \tilde\Omega$ as
an abstract manifold that can be embedded to $\p \Omega\subset \R^n$,
we have found the metric tensor on $
\p \tilde\Omega$ corresponding to the Euclidean metric of $\p \Omega$.
\smallskip
By the Cohn-Vossen rigidity theorem, intrinsically
isometric $C^2$-smooth surfaces that are boundaries of a strictly convex body
 are congruent in a rigid motion. For uniqueness,  
see e.g. \cite[Thm. V and VI]{Sacks} 
and also \cite{Hsu,Iaia}. This means
that the boundary data determines uniquely the map $T\circ f_m^{-1}$,
where $T$ is a rigid motion. Hence we can find the surface $T(\p \Omega)$ and on it
the map $T_*\Lambda_{\tilde \gamma}=
 T_*\Lambda_{\gamma}$.
Using the uniqueness of
of the isotropic inverse problem \cite{SU,Nachman},
we see that the boundary data determines $\gamma\circ T^{-1}$.
\hfill$\Box$ 
\bigskip

Note that the construction of the surface $\p \Omega\subset \R^3$ {\newtext 
from the intrinsic metric $h^E$} 
is a more delicate issue, see \cite{Nirenb,Other1}, hence we take care to avoid it.
\medskip

\section{A reconstruction algorithm and the  use of conformal flatness.}
In this section we consider the case $n=3$, even though the considerations
could be generalized for $n\geq 4$ by changing the 
Cotton-York 
tensor to Weyl tensor in our considerations (see Appendix).
As noted before, an actual construction of the isometric
embedding of an abstract manifold to Euclidean space is
complicated and thus we try to avoid it.
\smallskip

We want to find an anisotropic conductivity 
$\eta$ such that $R_{\tilde z,\eta}=\tilde R$
assuming that $\tilde R=(f_m)_*R_{z,\gamma}$
where $\gamma$ is an isotropic conductivity.
Clearly, {\newtext when $F_m:\Omega\to \tilde \Omega$ is 
diffeomorphism satisfying $F_m|_{\p \Omega}=f_m$,
the anisotropic conductivity $(F_m)_*\gamma$ is a solution
of the inverse problem, 
but it is not unique. }
However, we also know that $(F_m)_*\gamma$ 
has a conformally flat structure and this fact will help in solving
the inverse problem as we will see. 
Note that in principle,
one could start to solve the inverse problem by minimizing
over {\newtext  all pairs $(\Omega,\sigma)$ of smooth domains $\Omega\subset
\R^n$ and all isotropic conductivities $\sigma$ in $\Omega$.}
However, the minimization over domains is complicated, and our objective
is to find a reasonably simple minimization algorithm
where we minimize over conductivities in {\newtext the fixed model domain
$\tilde \Omega$} with an appropriately chosen cost function.
\smallskip

Let $\eta=(F_m)_*\gamma$ be a possibly anisotropic conductivity
in $\tilde\Omega$ such that $\gamma$ is isotropic.
As already noted 
it defines a Riemannian metric $g$ on $\tilde\Omega$, given by
\ba
[g_{jk}]_{j,k=1}^n=([g^{jk}]_{j,k=1}^n)^{-1},\quad g^{jk}= \det(\eta)^{1/(n-2)}\eta^{jk}
\ea
From now on  we will use the Einstein summation convention and omit the summation
symbols. 
As $F_m^{-1}:\tilde \Omega\to \Omega$ can be considered as coordinates, we see that in proper
coordinates the metric $g$ is a scalar function times
Euclidean metric, that is, $g$ is conformally flat.
This means that
\ba
g_{ij}(x)=e^{-2\sigma(x)}\overline g_{ij}(x)
\ea
where $\overline g_{ij}(x)$ is a metric with zero curvature
tensor (i.e.\ flat metric) and $\sigma(x)\in \R$.
By \cite{Eisen} (for original work, see \cite{Cotton}),
 the conformal flatness of the metric $g$ {\newtext in three dimensions} 
is equivalent to 
the vanishing of the Cotton-York tensor $\hbox{C}=\hbox{C}_{ij}$ corresponding to $g$
(see Appendix). 
Note that we can choose $\sigma=\frac {1}{2-n}\log \gamma$
and $\overline g=(F_m)_*(\delta_{ij})$. 
By \cite[formulae (28.18) and (14.1)]{Eisen}, $\sigma$ satisfies
a differential equation (with $n=3$)
\beq\label{eq: sigmalle}
\sigma_{ij}=
-\frac 1{n-2}\,\hbox{Ric}_{ij}+
\frac 1{2(n-1)(n-2)}\,g_{ij}\hbox{R}-
\frac 12\, g_{ij}g^{lm}\sigma_l\sigma_k,\quad i,j=1,\dots,n
\eeq
where $\hbox{Ric}_{ij}$ is the Ricci curvature tensor of $g$,
$R$ is the scalar curvature or $g$, and 
\ba
\sigma_k=\frac {\p \sigma}{\p x^k},\quad
\sigma_{ij}=\nabla_{e_i}\sigma_j-\sigma_i\sigma_j,\ \ \mbox{where $\quad e_i=\frac \p{\p x^i}$},
\ea
where $\nabla_{e_i}$ is the   covariant
derivative with respect to metric $g$.
Thus if $g$ is given, (\ref{eq: sigmalle}) is a second order
nonlinear differential equation for $\sigma$. 
By \cite[p. 92]{Eisen}, the equations
(\ref{eq: sigmalle}) satisfy the sufficient integrability conditions
to be locally solvable if and only if the Cotton-York tensor
vanishes. Note that the existence of the isotropic conductivity
$\gamma$ in $\Omega$ gives a solution for these equations.
\smallskip

\noindent
Consider now the following algorithm:
\medskip

\noindent
{\bf Data:} Assume that we are given $\p \Omega _m$,
$\tilde R=(f_m)_* R_{\gamma,z}$ and $z\circ f_m^{-1}$
on $\p \Omega _m$.
\medskip

\noindent
{\bf Aim:} We look for a metric $\tilde g$ corresponding to the conductivity $\tilde \gamma$
and $\tilde z$ such that on $\p \Omega _m$
$\tilde R= R_{\tilde \gamma,\tilde z}$ and
$\tilde z=(f_m)_*z$. 
\medskip

\noindent
{\bf Algorithm:}
\begin{enumerate}

\item Determine the two leading terms in the symbolic expansion of $\tilde R$.
They determine a contact impedance
 $\hat z$ and a metric $\hat g$ on $\p \tilde\Omega$
such that
if  $\tilde R= R_{\tilde \gamma,\tilde z}$
then
 $\tilde z=\hat z$ and $\tilde i^*(\tilde g)=\hat g$.

\item Form the ratio of reconstructed i.e. $\hat z$, and measured contact impedances
\ba
\hat r(\tilde x):=\frac 
{z(f_m^{-1}(\tilde x))}{\hat z(\tilde x)} ,\quad \tilde x\in \p \tilde\Omega.
\ea
Note that then 
\ba
\hat r(\tilde x)(f_m)_*(dS_E) = d\tilde S_E
\ea
since the contact impedances transformed as densities.

\item Let $dS_{\hat g}$ be the volume form of $\hat g$ on $\p \tilde\Omega$. Then
\ba
dS_{\hat g} = (\det \hat g)^{1/2}\, d\tilde S _E.
\ea
Define
\ba
\hat \gamma= (\det \hat g)^{1/2}\, \hat r.
\ea
With this choice
 $\hat \gamma$ will satisfy
$
\hat \gamma(\tilde x)=\gamma(f_m^{-1}(\tilde x))$ for $ \tilde x\in \p \tilde
\Omega.
$

\item Define the boundary value $\hat \sigma$ for the function $\sigma$
by
\ba
\hat \sigma(\tilde x)=  \frac {1}{2-n}
 \log \left(\hat \gamma (\tilde x)\right),
\quad \tilde x\in \p \tilde \Omega.
\ea

\item
Solve the  minimization problem
\ba
\min F_{\tau}(\tilde z,\tilde \sigma, \tilde \gamma)+\alpha\,
 H(\tilde z,\tilde \sigma,\tilde \gamma)
\ea
where $H(\tilde z, \tilde \gamma)$ is a regularization functional,
say
\ba
H(\tilde z,\tilde \sigma,\tilde \gamma)=\|\tilde z\|_{H^8(\tilde\Omega)}
+\|\tilde \gamma\|_{H^8(\tilde\Omega)}^2+\|\tilde \sigma\|_{H^8(\tilde\Omega)}^2,
\ea
$\alpha\geq 0$ is a regularization parameter, and
\ba
F_{\tau}(\tilde z,\tilde\sigma, \tilde \gamma)&=&
\|\tilde R-R_{\tilde \gamma,\tilde z}\|_{L(H^{-1/2}(\p \tilde \Omega))}^2
\\
& &+
\|\frac {\tilde z(\tilde x)}
{z(f_m^{-1}(\tilde x))}-\hat r(\tilde x)\|_{L^2(\p \tilde \Omega)}
\\
& &+
\tau\|\hbox{C}\|_{L^2(\tilde \Omega)}^2\\
& &+
\|\tilde\sigma|_{\p \tilde \Omega}-\hat \sigma\|^2_{L^2(\p \tilde \Omega)}
\\
& &+\sum_{i,j=1}^n \|\tilde\sigma_{ij}-
\left(-\frac 1{n-2}\hbox{Ric}_{ij}+
\frac 1{2(n-1)(n-2)}\tilde g_{ij}\hbox{R}-
\frac 12 \tilde g_{ij}\tilde g^{lm}\sigma_l\sigma_k\right)
\|_{L^2(\tilde \Omega)}^2
\ea
where $\tau\geq 0$, 
$\tilde g$ is the metric tensor corresponding to $\tilde \gamma$, 
$\hbox{C}=\hbox{C}_{ij}$ is the Cotton-York tensor of $\tilde g$,
and finally Ric  and $R$ are the Ricci curvature and
scalar curvature tensors  respectively of $\tilde g$.
\end{enumerate}

Note that above value of the Cotton-York tensor at $x\in \Omega$,
 $\hbox{C}_{ij}(x)$, the Ricci curvature tensors $R_{ij}(x)$, and 
the scalar curvature $R(x)$ depend on the values of the conductivity $\eta$
and its three first derivatives at $x$.

\begin{prop}\label{prop1}
Let $\Omega\subset \R^3$ be a bounded, strictly convex, 
$C^{\infty}$--domain. Assume 
that $\gamma \in C^{\infty}(\overline
\Omega)$ is an isotropic conductivity, $z\in C^\infty(\p \Omega)$,
$z>0$ be a contact impedance and $R_{\gamma,z}$ be the
corresponding Robin-to-Neumann map. 
Let
$\tilde\Omega$ be a model of the domain satisfying the same regularity 
assumptions as $\Omega$, and 
$f_m :\p\Omega\to \p\tilde\Omega$ be  a $C^{\infty}$--smooth diffeomorphism.

Assume that we are given $\p \tilde\Omega$, the values 
of the contact impedance $z(f_m^{-1}(\tilde x))$,
$\tilde x\in \p\tilde\Omega$,  and the map $\tilde R=(f_m)_*R_{\gamma,z}$.

Let $\tau\geq 0$. Then minima of $F_{\tau}(\tilde z,\tilde\sigma,\tilde \gamma)$ is zero and any minimizers $\tilde z$, $\tilde\sigma$ and $\tilde \gamma$
of
$F_{\tau}(\tilde z,\tilde\sigma,\tilde \gamma)$ satisfy $\tilde z=(f_m)_*z$
and there is a diffeomorphism $\tilde F:\Omega\to \tilde \Omega$ such that
$\tilde F|_{\p \Omega}=f_m$, $\tilde \gamma=\tilde F_*\gamma$ and $\tilde\sigma = -\log \tilde\gamma$.
\end{prop}

{\bf Proof.}  Assume first that $\tau>0$. 
The minimizer exists because of existence of $\Omega,\gamma,z$ and $\sigma$,
and the minimum is zero.
Let $\tilde z$, $\tilde\sigma$ and $\tilde g$ be some minimizers of $F_\tau$. 
{\newtext As then the Cotton-York tensor is zero and the equations
(\ref{eq: sigmalle}) are valid, if follows from \cite{Eisen},
that the metric}
$\overline g_{ij}=\exp(2\sigma(\tilde x))g_{ij}(\tilde x)$,
$x\in \tilde \Omega$ is flat. Since $R_{\tilde z,\tilde \gamma}=\tilde R$,
we  have $\tilde z=(f_m)_*z$ and 
the metric $\tilde g$ corresponding to $\tilde \gamma$
has to satisfy on the boundary  $i^*\tilde g=\hat g$.
This, and the vanishing
of $F_\tau$ imply that
 \ba
i^*\overline g&=&\exp(2\hat \sigma)i^* \tilde g=
\exp(2\hat \sigma) \hat g   
=\exp(2\hat \sigma)(f_m)_*(\gamma h_E)= \\
&=&
\exp(2\hat \sigma)\,\hat \gamma \,(f_m)_*(h_E)
=(f_m)_*(h_E).
\ea 
Consider now  $(\tilde \Omega,\overline g)$ as a Riemannian manifold.
As $\overline g$ is flat, we know that  $(\tilde \Omega,\overline g)$
can be embedded isometrically to domain $\Omega_0\subset \R ^n$.
Let $k:\tilde \Omega\to \Omega_0$ be this embedding.
Since $i^*\overline g=(f_m)_*(h_E)$, it follows from
the Cohn-Vossen rigidity theorem that the boundary
$\p \Omega_0$ and $\p \Omega$ are congruent in a rigid motion $T$
and $k\circ f_m=T|_{\p \Omega}$. {\newtext Then 
$(T^{-1}\circ k)_*\tilde \gamma$ is isotropic conductivity, the contact
impedances of $(T^{-1}\circ k)_*\tilde z$ and $z$ coincide, and
the Robin-to-Neumann maps of $(T^{-1}\circ k)_*\tilde \sigma$ and $\sigma$ coincide.
By the uniqueness of the isotropic inverse conductivity
problem \cite{SU}, $(T^{-1}\circ k)_*\tilde \gamma=\gamma$.}
This proves the claim in the case $\tau>0$.

Next, consider the case $\tau=0$. {\newtext Again, 
minimizer exists because of existence of $\Omega,\gamma,z$ and $\sigma$,
and the minimum is zero. Let $\tilde z$, $\tilde\sigma$ and $\tilde g$ be some minimizers. As the minima of $F_\tau$ is zero, the equations
(\ref{eq: sigmalle}) are valid.
By  \cite[p. 92]{Eisen}, the solutions $\sigma$ 
satisfy the integrability
conditions
\beq\label{eq: integrab}
\nabla_k\sigma_{ij}-\nabla_j\sigma_{ik}=\sigma_l \hbox{R}^l_{ijk},
\quad i,j,k=1,\dots,n
\eeq
that imply that the conformal covariant satisfies $R_{ijk}$ vanishes. 
Thus the Cotton-York tensor $\hbox{C}_{ij}$ 
is zero. This means that
the minimizers $\tilde z$, $\tilde\sigma$ and $\tilde g$ 
of $F_\tau$ with $\tau=0$
are also minimizers 
of $F_\tau$ with any $\tau>0$.}
\hfill$\Box$
\smallskip

\noindent One can think of $\tau$ as a regularization parameter: in general the solvability properties of equations (\ref{eq: sigmalle}) are sensitive to the compatibility conditions, i.e. the vanishing of the Cotton--York (or the Weyl tensor in higher dimensions).

To find the domain $\Omega$, we can continue the  above algorithm 
by applying the fact that conformally Euclidean manifold of dimeision $n$ 
can be conformally embedded to $\R^n$ 
in a  constructive way (cf.\ \cite{KKL}). 

\begin{enumerate}

\item[6.] {\newtext In the previous steps 1.--5.\ we have found
metric tensors $\tilde g$ and $\overline g=e^{2\tilde \sigma}\tilde g$ on $\tilde \Omega$
such that $\tilde g=F_*(g) $ and $\overline g=F_*(g^E)$ where
$g$ is the metric corresponding to
the metric $\gamma$ on $\Omega$, $g^E$ is the Euclidean metric 
on $\Omega$, and $F:\Omega\to \tilde \Omega$ is some diffeomorphism.

Let $y\in \tilde \Omega$ and find geodesics 
$\overline \mu_{y,\xi}(s)$ starting from $y$   
with respect to the metric $\overline g$. We parametrize these geodesics in such a way 
that 
 $\overline \mu_{y,\xi}(0)=y$    
and  $\p_s\overline \mu_{y,\xi}(0)=\xi$ is a unit tangent vector
of the tangent space $(T_{y}\tilde \Omega,\overline g)$. These geodesics correspond
to the halflines in $\R^3$ starting from some point $y_0\in \Omega$.
Let $J:(T_{y}\tilde \Omega,\overline g)\to \R^3$ be a linear 
isometry and define a map $\kappa:\tilde \Omega\to \R^3$ by setting 
\ba
\kappa(\overline \mu_{y,\xi}(s))\,=\,s\, J\xi,\quad s\geq 0.
\ea
Then $\kappa\circ F:\Omega\to \R^3$ is an affine isometry 
that extends to a rigid motion $T:\R^3\to \R^3$
with $T(y)=0$.
Thus
we can find $\kappa(\tilde \Omega)=T(\Omega)$,
$\kappa_*(\tilde \sigma)=T_*\gamma$, and
$\kappa_*(\tilde z)=T_*z$}
\end{enumerate}
Thus we have shown the followig reconstruction result.

\begin{coro}
{\newtext Let $\Omega,\gamma,z,\tilde \Omega$ and $f_m$ be
as in Proposition \ref{prop1}. Assume that we are 
given $\p \tilde\Omega$, the contact impedance $z(f_m^{-1}(\tilde x))$,
$\tilde x\in \p\tilde\Omega$,  and the Robin-to-Neumann map $\tilde R=(f_m)_*R_{\gamma,z}$.
Then the algorithm 1.--6.\ determines $\Omega$, 
$\gamma$, and $z$ up to a rigid motion $T:\R^3\to \R^3$.}
\end{coro}

We intend to investigate the numerical implementation of the method and give numerical test results
in part II of this paper.



\section*{Appendix.} Here we define the conformal curvature tensors.
{\newtext We say that a metric $g_{ij}$ in a domain $\Omega\subset
\R^n$ is conformally
flat if there is a scalar function $a(x)>0$ such that
the curvature of tensor of $a(x)g_{ij}(x)$ is identically zero.

First, let $\gamma$ be an isotropic conductivity, i.e., a 
smooth positive function in $\Omega$ and $F:\Omega\to \tilde \Omega$
be a diffeomorphism.}
Let $\eta=F_*\gamma$ be a possibly anisotropic conductivity
in $\tilde\Omega$. 
It defines a Riemannian metric $\tilde g$ on $\tilde \Omega$, given by
\ba
[\tilde g_{jk}]_{j,k=1}^n=([\tilde g^{jk}]_{j,k=1}^n)^{-1},\quad \tilde g^{jk}= \det(\eta)^{1/(n-2)}\eta^{jk}
\ea
As $F^{-1}:\tilde \Omega\to \Omega$ can be considered as coordinates, we see that in proper
coordinates the metric $\tilde g$ is a scalar function times
Euclidean metric, that is, $\tilde g$ is conformally flat.

{\newtext Next we consider a general metric tensor $g_{ij}$ and
recall facts concerning its conformal flatness}
Note that below we use the Einstein summation convention and omit the summation
symbols when possible. The following tensors are related to conformal
flatness.

\begin{itemize}

\item[(a)] Assume that $n=3$. {\newtext Then the conformal covariant, given 
in terms of curvature tensors (see the explanation on notations
below) is
\ba
R_{ijk}=
\nabla_kR_{ij}-\nabla_jR_{ik}+\frac 1{2(n-1)}(g_{ik}\nabla_jR-
g_{ij}\nabla_kR).
\ea
{\newtext In the three dimensional case $R_{ijk}$
defines a tensor that can be considered as vector valued
2-form $R_{ijk}dx^j\wedge dx^k$. Operating with the Hodge operator $*$ 
to this 2-form, we obtain the Cotton-York tensor,}
\ba
\hbox {C}_{ij }
=g^{k p } g^{l q} \nabla_k(R_{li}-\frac 14 R\,g_{li})\epsilon_{pq j},
\ea
where $\epsilon_{pqj}$ is the Levi-Civita permutation symbol.}

\item[(b)] Assume that  $n\geq 4$. Then the Weyl  tensor is
\ba
W_{i j kl}&=&R_{i j kl} +
\frac 1{n-2}(g_{i l}R_{kj }+g_{j k}R_{kj }-g_{i k}R_{lj }-g_{j l}R_{ki} )+
\\ & &+
\frac 1{(n-1)(n-2)}(g_{i k}g_{lj }-g_{i l}g_{kj })R.
\ea
\end{itemize}
{\newtext The crucial fact related to our considerations is that 
the metric $g$ is conformally flat
if and only if in the dimension $n=3$ the Cotton-York tensor vanishes
and  in the dimension $n=4$ the Weyl tensor vanishes,
see \cite[p. 92]{Eisen} or \cite{Cotton,Thomas,Aldershey}}.

Above, $R_{i j kl}$ is the Riemannian curvature tensor,
\ba
R_{i j kl}=\frac {\p }{\p x^k }
 \Gamma^i _{j l}-
\frac {\p }{\p x^l } \Gamma^i _{j k}+ \Gamma^p_{j l} 
\Gamma^i _{pk} - \Gamma^p_{j k} 
\Gamma^i _{pl},\quad R^p_{j kl}=g^{pi} R_{ijkl}
\ea
where $\Gamma^i _{j k}$ are Christoffel symbols,
\ba
\Gamma^i _{j k}=\frac 12 g^{pi }(\frac {\p g_{j  p}}{\p x^k}
+\frac {\p g_{k p}}{\p x^j }-\frac {\p g_{j  k }}{\p x^p }),
\ea
$R_{i j }$ is the Ricci curvature tensor, 
$
R_{i j }=R^k_{i jk },
$
and $R$ is the scalar curvature 
$
R=g^{i j }R_{i j }.
$
Finally, above $\nabla_k$ is the covariant derivative that is
defined for a (0,2)-tensor $A_{il}$  and a (0,1)-tensor $B_{l}$ by
\ba
\nabla_kA_{li}=\frac{\p}{\p x^k}A_{li} -\Gamma^p_{kl}   A_{pi}
-\Gamma^p_{ki}A_{lp},\quad\quad
\nabla_kB_{l}=\frac{\p}{\p x^k}B_{l} -\Gamma^p_{kl}   B_{p}.
\ea

\noindent
{\bf Acknowledgements.} {\newtext The research was supported by
Finnish Centre of Excellence in Inverse Problems Research  
(Academy of Finland CoE--project 213476).}


\begin{thebibliography}{99}


\bibitem{adler96b}
Adler A, Guardo R,  Berthiaume Y.
Impedance imaging of lung ventilation: {D}o we need to account for
  chest expansion?,
 {\em IEEE Trans. Biomedical Engineering}, {\bf 43} (1996), 414--420.

\bibitem{Aldershey}
Aldersley S. Comments on certain divergence-free tensor densities in a 3-space, {\em J. Math. Phys.} {\bf 9}(1979), vol. 20, 1905-1907. 

\bibitem{alp} Astala, K.,  Lassas  M.,  P\"aiv\"arinta L.   
Calder\'on's inverse problem for anistropic conductivity
in plane, 
{\em Comm. Partial Differential Equations}  {\bf 30}  (2005), 207--224.

\bibitem{ap} Astala, K.  P\"aiv\"arinta  L. Calderon's inverse conductivity problem in the plane,  {\em Ann. of Math.} (2)  {\bf 163}  (2006), 265--299.


\bibitem{BU}
Brown, R., Uhlmann, G.
Uniqueness in the inverse conductivity problem for nonsmooth
conductivities in two dimensions.
{\em Comm. Partial Differential Equations}  {\bf 22} (1997), 1009--1027.



\bibitem{ca} Calder\'on, A.  On an inverse boundary value problem, Seminar on Numerical Analysis
and its Applications to Continuum Physics (Rio de Janeiro, 1980), {\em Soc. Brasil Mat.} (1980), 65--73, Rio de Janeiro

\bibitem{chengea89}
Cheng, K.-S., Isaacson  D., Newell  J.C., Gisser D.G.  Electrode models for electric current computed tomography, {\em IEEE Trans. Biomed. Engrg.} 3, {\bf 36}
(1989), 918--924.






\bibitem{Cotton} 
 Cotton E. Sur les varietes a trois dimensions,
{\em  Ann. Fac. Sci. Toulouse} 
 (1899),  no. 4, 385--438.

\bibitem{Eisen} 
 Eisenhart L. {\em Riemannian Geometry}, Princeton University Press,
Princeton, N.J. 1925, 1977, 



\bibitem{GLU}
Greenleaf, A., Lassas, M., Uhlmann,  G. 
The Calderon problem for conormal potentials, I: Global uniqueness and reconstruction,  {\em Comm. Pure Appl.} Math  {\bf 56} (2003), 328-352.


\bibitem{Hsu}
 Hs\"u, C.  Generalization of Cohn-Vossen's theorem.  
{\em Proc. Amer. Math. Soc.}  {\bf 11}  1960 845--846.


\bibitem{Iaia}
Iaia, J.
Isometric embeddings of surfaces with nonnegative curvature in $R\sp 3$.
{\em Duke Math. J.} {\bf 67} (1992), 423--459

\bibitem{Siltanen3}
Isaacson D., Mueller J., Newell J., Siltanen S.
Reconstructions of chest phantoms by the d-bar method for electrical impedance tomography,
{\em IEEE Transactions on Medical Imaging} {\bf 23}(2004), 821--828.



\bibitem{Jost}
Jost, J. {\em Riemannian Geometry and Geometric Analysis}, Springer, 1995; 
2001


\bibitem{gersing96}
Gersing E., Hoffman B.,  Osypka M.
\newblock Influence of changing peripheral geometry on electrical impedance
  tomography measurements.
\newblock {\em Medical \& Biological Engineering \& Computing},{\bf 34} 
(1996), 359--361.


\bibitem{kolehmainen97e}
Kolehmainen V, Vauhkonen M, Karjalainen P.,  Kaipio J.
Assessment of errors in static electrical impedance tomography with
  adjacent and trigonometric current patterns.
 {\em Physiological Measurement}, {\bf 18}(1997), 289--303.


\bibitem{kaipio00a}
Kaipio J., Kolehmainen V., Somersalo E., Vauhkonen M. 
 Statistical inversion and {{M}onte} {{C}arlo} sampling methods in
  electrical impedance tomography.
 {\em Inverse Problems}, {16}(2000), 1487--1522.



\bibitem{Kolehmainen}
Kolehmainen V. , Lassas, M., Ola P.  Inverse conductivity problem with an imperfectly known boundary, {\em SIAM J. Appl. Math.}
{\bf 66}  (2005), 365--383. 

\bibitem{KKL}
Katchalov, A.  Kurylev Y., Lassas M. {\em
 Inverse boundary spectral problems}. Chapman Hall/CRC Monographs and Surveys in Pure and Applied Mathematics, {123}., 2001. xx+290 pp

\bibitem{LU}
 Lassas, M., Uhlmann, G. On determining a Riemannian manifold from the Dirichlet-to-Neumann map.  {\em Ann. Sci. Ecole Norm. Sup.}  {\bf 34}  (2001), 771--787.

\bibitem{LUT}
 Lassas, M., Taylor, M.E., Uhlmann, G. The Dirichlet-to-Neumann map for 
complete Riemannian manifolds with boundary.  {\em Comm. Anal. Geom.} {\bf  11}  (2003), 207--221.

\bibitem{LeeUhlmann}
Lee, J., Uhlmann, G.
Determining anisotropic real-analytic conductivities by boundary measurements.
{\em Comm. Pure Appl. Math.} {\bf 42} (1989), 1097--1112.


\bibitem{Lionheart}
Lionheart, W. Boundary shape and electrical impedance tomography.  
{\em Inverse Problems}  {\bf 14}  (1998), 139--147.



\bibitem{Nachman}
 Nachman, A. Reconstructions from boundary measurements. 
{\em Ann. of Math.} (2)  {\bf 128}  (1988), 531--576.



\bibitem{na2} Nachman, A.  Global uniqueness for a two--dimensional inverse boundary value problem,
{\em Ann. Math.} {\bf 143} (1996), 71--96.

\bibitem{Nirenb}
 Nirenberg, L. The Weyl and Minkowski problems in differential geometry in the large.  {\em Comm. Pure Appl. Math. } {\bf 6} (1953), 337--394.



\bibitem{PPU} P\"aiv\"arinta L., Panchenko A., Uhlmann G. 
Complex geometrical optics solutions for Lipschitz conductivities
Rev. Mat. Iberoamericana.  {\bf 19}  (2003), 57--72.


\bibitem{Pekonen}
 Pekonen, O. On the variational characterization on conformally flat $3$-manifolds.  {\em J. Geom. Phys.}  {\bf 7}  (1990), 109--117


\bibitem{Other1}
Pogorelov, A.
The rigidity of general convex surfaces. (Russian)
{\em Doklady Akad. Nauk SSSR} (N.S.) {\bf 79}, (1951) 739--742.


\bibitem{Sacks}
Sacksteder, R. The rigidity of hypersurfaces.  {\em J. Math. Mech.} {\bf 11}  1962, 929--939.


\bibitem{Siltanen2}
Siltanen S., Mueller J., Isaacson D.
An implementation of the reconstruction algorithm of A. Nachman for the 2-D inverse conductivity problem,
{\em Inverse Problems} {\bf 16}(200), 681-699. 

\bibitem{Siltanen1}
Mueller J., Siltanen S.
Direct reconstructions of conductivities from boundary measurements,
{\em SIAM Journal of Scientific Computation} {\bf 24}(2003), 1232--1266.



\bibitem{Somersalo92}
Somersalo E., Cheney M., Isaacson D.  
Existence and uniqueness for electrode models for electric current computed tomography.
{\em SIAM J. Appl. Math.}, {\bf 52}(1992), 1023--1040.


\bibitem{sy} Sylvester, J.  An anisotropic inverse boundary value problem, {\em Comm. Pure Appl. Math.} {\bf 43} (1990),
201--232.


\bibitem{SU}
 Sylvester, J., Uhlmann, G. A global uniqueness theorem for an inverse boundary value problem.  {\em Ann. of Math.} (2)  {\bf 125}  (1987), 153--169.



\bibitem{Thomas} 
Thomas J. Conformal Invariants
{\em Proc. Natl. Acad. Sci. USA.} 1926 12(6), 389-393.



\bibitem{U}
Uhlmann G.: Inverse boundary value problems for partial
differential equations. 
{\em Proceedings
of the ICM.} Vol. III (Berlin,
1998).
Doc. Math. Vol. III, 77--86.




\end{thebibliography}
\end{document}